\documentclass[pamq]{ipart}

\startlocaldefs

\theoremstyle{definition}

\theoremstyle{remark}

\numberwithin{equation}{section}

\def\IC{\mathbb{C}}
\def\IH{\mathbb{H}}

\def\IR{\mathbb{R}}

\def\IO{\mathbb{O}}

\def\CM{\mathcal{M}}
\def\CG{\mathcal{G}}

\def\cm{{\,\mbox{cm}}}
\def\TeV{{\,\mbox{TeV}}}

\endlocaldefs

\pubyear{2019}
\volume{0}
\issue{0}
\firstpage{1}
\lastpage{1}
\begin{document}

\begin{frontmatter}

\title{From Algebraic Geometry to Machine Learning}

\begin{aug}

\author{\fnms{Michael R.} \snm{Douglas}\ead[label=e1]{mdouglas@scgp.stonybrook.edu}}
\address{Simons Center for Geometry and Physics,
Stony Brook University, Stony Brook NY 11794\\\printead{e1}}
\end{aug}





\end{frontmatter}

\begin{abstract}
David Mumford made groundbreaking contributions in many fields, including
the pure mathematics of algebraic geometry and the applied mathematics of machine learning
and artificial intelligence.  His work in both fields influenced my career at several key moments.
\end{abstract}

\maketitle

David Mumford had a significant influence on my career.  While many mathematicians 
and mathematical scientists can say that,
in my case it requires explanation, as I have only met him once, and briefly at that.
I believe I have never referenced his work before now, nor was that an evident omission;
his influence on me was less through specific results and more through the clarity and impeccable taste of his work.
Indeed, when I was invited to contribute to this special issue, my first thought was,
how did the organizers know about this?  Eventually I asked them, but let me save their reply for the end.

Let me begin with the story of the first and so far only time I met David Mumford.
It was the spring of 1983, and I was a senior at Harvard, about to graduate with a physics degree.
I was deciding between Caltech and Princeton for graduate school, and had just returned from
a visit to Caltech.  Although I was fascinated by fundamental physics and its connections with
mathematics, I was equally fascinated by computer science, cognitive science and especially artificial intelligence.
And that year I was mostly taking classes related to the latter direction, on topics such as analytic
philosophy and the psychology of memory, while reading everything I could find about AI.

During my visit to Caltech I attended a great class.  Richard Feynman, John Hopfield and Carver Mead
had gotten together to teach a course on the physics of computation.  This course, and the three courses
that the three professors taught individually the next year, was arguably one of the seminal events in
the history of computer science.  Feynman gave the first series of lectures on his new concept of a quantum
computer.  Mead introduced a new approach to integrated circuit design which combined digital and analog concepts, 
an approach often referred to now as neuromorphic computing.  And Hopfield gave the first series of lectures
on his new model of a neural network, based on ideas from the theory of spin glasses.

Now at that moment, the future of fundamental physics was not looking so bright.
The proton-antiproton collider at Cern was running, but after the discovery of the W and Z bosons there was no clear evidence for the Higgs boson or for any other particles beyond the Standard Model (having found the Higgs boson in 2012, we are back in this situation, by the way).  The remarkable discoveries that would reshape cosmology such as dark energy and the detailed measurement of the cosmic microwave background were still far in the future.  On the theoretical side, the burst of progress on Yang-Mills theory that had fueled the development of the Standard Model had largely stalled, as none of the approaches to nonperturbative calculations had yet borne fruit.  Supersymmetry was considered a wild speculation, while few people had even heard of superstring theory.  While I did not know much of this at the time, I was well advised and soon decided that I would start my graduate studies in the interface between physics and computation, leaving open the possibility of going into fundamental physics should the prospects look brighter.

So, I came back from Caltech with papers and notes from this fascinating class, including Hopfield's original preprint on his model \cite{Hopfield}.
And for some reason I don't quite remember but which I think was administrative, I had to visit the math department offices on an upper floor of the Science Center.  As I sat waiting for a signature on a form, I was reading Hopfield's paper, when I was interrupted by a distinguished looking professor asking me what I was reading.  When I told him, he asked if he could copy the paper, and I gave it to him.  He told me his name, David Mumford, and I thought little more of it.

At Caltech the next fall, I took Feynman's and Hopfield's courses, and prepared to work on neural networks.  I even wrote a paper for Feynman's class on quantum computing, though I was not farsighted or brave enough to pursue that further.
I also joined Gerry Sussman, one of the early pioneers of AI who was on a sabbatical visit from MIT, on his project to build a computer he called the Digital Orrery.  This was a special purpose computer designed to integrate the equations of celestial mechanics governing the Solar System for periods of up to a billion years, to explore questions of chaos and stability \cite{Orrery}.  With this computer we showed that the orbit of Pluto is chaotic \cite{SW}, and I briefly considered switching to astrophysics, but before getting very far I learned about a discovery which overshadowed celestial mechanics and even neural networks.  This was the anomaly cancellation in ten-dimensional superstring theory, discovered in the summer of 1984 by Michael Green and John Schwarz.  John was then a senior research fellow at Caltech and this news electrified the physics world.  Many graduate students in theoretical particle physics immediately began working on superstring theory, and I was swept up in the excitement, which turned out to be amply justified.

Back then, one of the central problems of superstring theory was to show that the perturbative expansion was well defined and finite at all orders.  Making concrete computations required doing integrals over the moduli space of Riemann surfaces, while showing the physical consistency of the theory required demonstrating factorization.  This meant looking at every possible limit in moduli space, appealing to mathematical results that showed that every singular limit corresponded to a degeneration in which a surface breaks up into a pair of lower genus surfaces with punctures,
and showing that the singular behavior of the integral had an interpretation in terms of the propagation of strings between the pair of surfaces.  The most attractive formulation of this moduli space, which both exhibited its holomorphic structure and made the limiting behaviors clear, was the Deligne-Mumford compactification.
And in doing concrete calculations, one was soon initiated into the wonders of automorphic functions and especially theta functions.  By my last year at Caltech I had my own copies of Mumford's Tata lectures \cite{MT}, and their combination of very explicit calculations placed in the service of a clear general picture was appealing to me.

I moved to Chicago in the fall of 1988, mainly to work with Dan Friedan and Steve Shenker.  I was inspired by their bold proposal \cite{FS} to make an all-orders definition of string theory which could be reinterpreted as a nonperturbative definition.  We discussed approaches to doing this, and it was here that I drew on the lessons that I learned from the Tata lectures as well as Mumford's lectures on curves \cite{MC}.  The most important lesson was the great difficulty of working with an explicit definition of the moduli space for genus three and greater.  This got us interested in alternate definitions of string theory which did not need such a definition.  Soon we were looking at simplifications of the string theory problem, especially Polyakov's noncritical string and the work on discrete two-dimensional quantum gravity and its connection to matrix models developed by Migdal with his students Kazakov and Kostov.

Our search for a solvable nonperturbative string theory led us to formulate and solve the double scaling limit 
of matrix models \cite{double}.   Its main relevance for my story
here was that through its reformulation in terms of the KdV hierarchy \cite{Douglas:1989dd}, it led to
Witten's conjecture \cite{Witten:1990hr} governing the intersection cohomology of moduli space, a question which had been formulated by Mumford in \cite{ME}.  As I slowly absorbed this development, I realized that I would do well to study Mumford's work more deeply.

Despite all the excitement and progress in string theory, I retained my fascination with cognitive science and AI.
This was not an easy interest to keep up, especially since there was zero overlap between it and superstring theory.
But there was a growing overlap between computer science and statistical physics, which I learned about during
visits to the \'Ecole Normale and to MIT in 1990.

When I returned, I heard that Mumford had changed fields, from algebraic geometry to statistics and machine learning.
By this time I had come to think of him as a lodestone, a trustworthy guide to mathematics which could be used to solve real problems and at the same time was of the highest conceptual significance.  So I was excited to see how such a mathematician would treat these topics.  I recall reading \cite{MN}, which I found stimulating but which, a bit to my surprise, contained no equations.  I was not ready to go that far, and instead studied phase transitions in combinatorial optimization, a topic which had been greatly advanced by Giorgio Parisi and Marc M\'ezard's group at the \'Ecole Normale,
and in which physics techniques have been very powerful. \cite{MM}

My interest in AI was again interrupted by revolutionary progress in superstring theory,
in this case the Seiberg-Witten solution for $N=2$ supersymmetric Yang-Mills theory.  Rutgers and Princeton
instantly became the center of what would soon be called the second superstring revolution.
Algebraic geometry and moduli spaces of Riemann surfaces again became central, but in a different role,
as moduli spaces of vacua for gauge theories and string theory compactifications.  Soon the Tata lectures were
back on my desk, and we were again looking at all possible singular limits in moduli spaces, though
now the physical interpretation was very different.  Let me cite \cite{AD} as one of many examples -- here
the degeneration limit of a genus two Riemann surface to a pair of elliptic curves had the physical interpretation
of an $N=2$ superconformal field theory, a novel and still somewhat mysterious beast.
But, important as this story is, let us skip ahead to the discovery of the Dirichlet brane or D-brane.

Joe Polchinski's paper \cite{Pol} appeared in 1995, and suddenly many of our questions about superstring duality and M theory had clear answers.  One of the first new observations was that, since the first step in defining a D-brane was to promote the fields parameterizing its position to matrices, one expected noncommutative geometry to be relevant.  I was very struck by this \cite{Douglas:1996vj} and it was a major influence leading me to spend the fall of 1997 at the IHES.  While the immediate product of my visit was a well known joint work with Alain Connes and Albert Schwarz \cite{Connes:1997cr}, I also had many conversations with Maxim Kontsevich whose somewhat different conception of noncommutative geometry was totally unknown in the string theory community.  On my return I continued studying the geometry of D-branes, 
following the physics progression of dealing successively with systems with less and less supersymmetry.
By 2000 we were focusing on the case of $N=1$ supersymmetry in four dimensions.

Good mathematical references on supersymmetry include \cite{Freed} and the book \cite{Clay}, which tells  
the story I am about to summarize. 
But let me give a brief taste of the physics here.
For many reasons, work on ``beyond the Standard Model physics'' has largely started from the hypothesis
that at distances much shorter than what we can now probe at accelerators (about $10^{-18} \cm$),
but much longer than the Planck length at which our concepts of space-time must change (about $10^{-33} \cm$),
physics is described by a quantum field theory with $N=1$ supersymmetry (or simply ``$N=1$ field theory'').
Such a theory corresponds to the following mathematical data (\cite{Clay} \S 5.4): a K\"ahler manifold $\CM$
with a holomorphic symplectic action of a group $\CG$, and an invariant holomorphic function\footnote{
In supergravity $W$ is a section of the line bundle whose curvature is the K\"ahler form.}
$W$ on $\CM$ called
the superpotential.  Its moduli space (of vacua) is then the space of solutions of $\nabla W=0$ in the
symplectic quotient $\CM//\CG$.  The restriction to solutions and the quotient commute, and in the physics literature
are referred to as ``solving the F and D-term conditions'' respectively.
These connections between physics and mathematics were established in the 80's \cite{Hitchin:1986ea}
and so the relations between symplectic quotient and GIT quotient and to Mumford's work were very well known to us.

Now, the theory of compactification of superstring theory on a Calabi-Yau manifold $M$ (summarized in \cite{Douglas:2015aga}),
tells us that the moduli space of the structures postulated in compactification --
Ricci flat metrics and hermitian Yang-Mills (HYM) connections on $M$, Dirichlet branes on $M$ (more about this shortly),
and other structures -- must be the moduli space of an $N=1$ field theory.  And one can make many formal connections
between math and physics at this level.  For example, the Donaldson-Uhlenbeck-Yau theorem, that a bundle $V$
on $M$ admits an irreducible HYM connection if and only if it is $\mu$-stable, can be related in physics terms
to D-flatness conditions which depend on the K\"ahler class of $M$. \cite{Sharpe}

Another entry point into the mathematics was the theory of quiver representations.  This entered almost from
the start with the relation between D-branes on orbifolds and the work of Kronheimer and Nakajima \cite{Douglas:1996sw},
and was a constant theme.  In studying quiver theories we learned about the work of King on $\theta$-stability \cite{King},
and with Mumford's work in mind we hypothesized that stability should be a central part of the discussion.
By combining these various ideas my students Bartomeu Fiol, Christian R\"omelsberger and I came up with a 
physics concept combining $\mu$-stability and $\theta$-stability which we called $\Pi$-stability,
after the standard notation $\Pi$ for the basis of central charges which entered the condition. \cite{Douglas:2000ah}

At this point Kontsevich's insights became a very fruitful guide.  One example which I had learned from him at IHES
was the interpretation of the holomorphic Chern-Simons action as a superpotential, and the conjecture that the
obstruction theory of any holomorphic object on a Calabi-Yau (say, for families of holomorphic curves) would be
governed by a superpotential.  He had also pointed out the relevance of derived categories and autoequivalences,
which was taken up by many physicists and mathematicians.  

The big question which Kontsevich's work raised was to understand homological mirror symmetry in physics terms.
Kontsevich's conjecture \cite{K} equated two categories, the derived category of coherent sheaves on $M$ or 
${\mbox DCoh}(M)$, and the Fukaya category on its mirror manifold $W$.  By this point we saw that these
were both categories of D-branes, and had many ideas about what physics could add to the discussion.  
Still, the heavy formalism of
the derived category was an obstacle.  Richard Thomas had written an excellent introduction \cite{Thomas:1999ic}
and I had many invaluable discussions with him and Paul Seidel about the topic.
One point I took from them was the difficulty of even defining the Fukaya category, which led me to focus on 
${\mbox DCoh}(M)$.  But the other was our shared belief in the importance of stability for the discussion.

I then went and asked many mathematicians how to define stability for objects in a derived category.  The
universal answer was that this was not possible; stability required a concept of subobject which did not exist there.
Rather, one should try to define a $t$-structure and find an abelian subcategory.  While sensible enough, the physics
did not suggest any way to do this.  Thus I went a different way, based closely on the physics.  While we could not
derive much from the 
general definition of D-branes, what was physically guaranteed to exist was an assignment of ``central charge''
(a complex character on the K theory)
and ``grading'' (a lifting of its phase to $\IR$)
for each stable object in the derived category, and a process of ``binding/decay'' by
which each distinguished triangle in the derived category could lead to a wall-crossing and change of stability
when three central charges align in the complex plane.
These ideas led to my proposal for $\Pi$-stability in the derived category
 \cite{Douglas:2000gi}, a work with Paul Aspinwall exploring the proposal physically
\cite{Aspinwall:2001dz}, and eventually my 2002 ICM lecture
\cite{Douglas:2002fj} which brought it to the attention of mathematicians,
particularly Tom Bridgeland who went on to turn it into rigorous mathematics \cite{Bridgeland}.

As I was accomplishing my goal of turning these problems over to the mathematicians, 
I started to get back into AI.  In addition to catching up with
the ``good old-fashioned AI'' school I had studied with Gerry Sussman, and catching up with the physics of combinatorial
optimization, I was finding that many of the pioneers of neural networks in the 90's were  
drawing increasingly on the theory of statistical inference \cite{MacKay}.  And at the 2002 ICM I heard David Mumford
himself expound on the subject \cite{MP}.  This talk may have been my first exposure to generative models,
which in our current age of data science have become extremely popular, but were almost unknown to theoretical physicists
at the time (a notable exception was Vijay Balasubramanian \cite{VB}). Not long after, 
Lee and Mumford developed the hypothesis that the architecture of the cerebral cortex, with its
multiple feedforward and feedback connections, was implementing a type of Bayesian statistical analysis
\cite{LM}.   This was also ahead of its time, and now it is one of the hot topics in cognitive science \cite{Clark}.

Although fascinating, there was a major obstacle to my working on AI: it still had almost no overlap
with string theory and my whole field of ``formal physics.''   And I was not in a position to follow Mumford's
example, to retire from string theory and devote myself to a second career.  For one thing, I was now the director
of the New High Energy Theory Center at Rutgers, working to rebuild and maintain the group after the departure
of two of its most prominent members.   But more importantly, we all looked forward to
a new golden age of particle physics when the Large Hadron Collider turned on, scheduled for 2008.
Either we would discover supersymmetry, or possibly something unexpected and even more exciting.
It was not at all the time to leave the field; rather it was time to prepare for contact between string theory and new data.

With these various motivations, I looked for a way to try to make predictions from string theory for the LHC
and other upcoming experiments, and to turn the problem into one of statistics,
leading to \cite{Douglas:2003um}.  Despite the way I came to it, 
to give a simple explanation of this I should start from neither string theory nor AI.  Rather, there was already a
field of ``quantum cosmology,'' in which the researchers extrapolated from better established concepts
such as inflation and the theory of phase transitions, to develop a framework to describe not just the Big Bang,
but whatever came before.  This is particularly interesting in a theory with many solutions, as one of the jobs
of cosmology is to try to explain which of the solutions are the preferred candidates to describe our universe.
This is usually phrased as the computation of the ``measure factor,'' an {\it a priori} probability distribution
over the solutions of the theory.

Now quantum cosmology had started well before string compactification and there were no strong internal
reasons to believe that the theory of quantum gravity should have many solutions, but there was an (in)famous
argument which required it to have many solutions, the anthropic solution to the cosmological constant problem.
As explained in \cite{Weinberg} (an update and list of references is in \cite{Polchinski:2006gy}),
in a theory of quantum gravity, the cosmological constant (the energy of empty space-time) gets
contributions from the quantum fluctuations of all of the fields in the theory, which are far larger in magnitude
than the very small observed value.  There are strong arguments that no physical mechanism can solve 
this problem directly, but it is also the case that if the universe had such a large
cosmological constant, it would either expand or contract so rapidly that there would not be enough time
for life to evolve to observe it.
Thus if we could argue that there were many solutions of quantum gravity with different values of the cosmological
constant, then as long as there are any solutions which are compatible with our existence, we could argue that the fact
that we observe such a solution requires no more explanation -- it is a selection effect.
This is a particular case of the anthropic principle by which one can go on to explain many other
constraints on the observed laws of physics, as required to get the nontrivial physics required for life.  
In this generality the principle is heartily disliked by physicists, who would much rather derive the
laws of physics from more fundamental principles.  Still, it is philosophically sound -- and given the strong
arguments that there is no other solution to a central problem, the cosmological constant problem,
one can find it compelling.

These general ideas had been around since the late 80's, but what brought them to the fore was the discovery in the 90's
that the expansion of the universe, discovered by Hubble in the 20's, was accelerating, in a way that could be
simply fit by postulating a small non-zero cosmological constant.  Before that, physicists had not given up on the
idea that some unknown mechanism was setting the cosmological constant to zero.  However, the idea of a mechanism
that sets it to a small non-zero constant value seemed even less plausible than the anthropic argument.  Then
in 2000, Bousso and Polchinski argued that string theory had the multiplicity of solutions required by the argument.
\cite{Bousso:2000xa}  Their original arguments were somewhat heuristic (even by physics standards), but combined with the
established fact that none of the other ingredients in string compactification were unique (there are many
Calabi-Yau manifolds, many holomorphic vector bundles on each, {\it etc.}) they convinced me that one could not
find the ``correct'' solution of string theory {\it a priori}; rather one needed to consider the ensemble of solutions
and treat the problem of finding the relevant ones and their likely predictions as a problem in statistical inference.

Without going into details, a major part of this problem is to enumerate solutions to $\nabla W=0$ as introduced
earlier, and a reasonable approach (well motivated by the Bousso-Polchinski work) was to consider an appropriate
random ensemble of superpotentials $W$, for example the period of a randomly drawn homology three-cycle,
and find the statistics of the critical points.  With this in mind I had started scouting out the mathematical literature,
and sometime in 2002 I ran into \cite{BSZ}, which studied zeroes of random sections of a line bundle, very similar
to my problem.  I was delighted to learn that Steve Zelditch would be at ICM 2002, and that is where our
collaboration (with Bernie Shiffman as well) began, leading to \cite{DSZ} and many subsequent works.

These ideas were exciting and successful enough to more or less supplant my original motivation of
learning statistics to work on AI,
at least as long as the prospect of new LHC physics was in sight.  However, in a story I have told elsewhere 
\cite{Douglas:2012bu},
following them up led to serious criticisms of the widely accepted arguments that supersymmetry would be
discovered at LHC -- not that one could prove it would not, but I came to believe that string theory did not
inevitably predict supersymmetry at LHC \cite{Douglas:2006es}, and that
one should give equal weight to the possibility that it would not be discovered.  Thus, as data started coming in and
superpartners were not found, I concluded in 2011 that this would be the likely outcome.  
I should stress that we do not know the final outcome yet --
the LHC has many more years to run, and more discovery potential -- but the original predictions from
``natural supersymmetry'' have been falsified, and (at least in my opinion  \cite{Douglas:2012bu}) 
the arguments that replace them suggest that supersymmetry should be discovered at energies between
30 and $100 \TeV$.
Such experiments are possible, but the larger accelerators they require would not come online until 2040 at the 
earliest.

Following a combination of these intellectual motivations and
personal reasons, in 2012 I left physics to pursue my interest in statistical inference in a very different way,
by joining Renaissance Technologies, one of the first and most successful quantitative hedge funds.
Even this step did not take me away from the beneficial influence of David Mumford, and I found his textbook 
with Desolneux \cite{MP2} very helpful as I came to terms with the vast practitioner's literature on machine
learning and finance.  But it does take me beyond the scope of this homage.  

Finally, I promised to say why the editors asked me to contribute, so I asked them.  The answer was that they felt
that stability conditions are an example of where ideas from string theory have come to influence most profoundly the areas of math that Mumford cared about and worked on,  and they wanted a contribution in the original spirit of the math physics.
So I'm delighted to have carried on his interests in this way, and for this chance to thank him for his broader contributions as well.

\bigskip

\bibliographystyle{amsplain}

\begin{thebibliography}{10}

\bibitem{Orrery}  ``A Digital Orrery,'' 
J. H. Applegate, M. R. Douglas {\it et al}, 
IEEE Transactions on Computers, {\bf C-34}, pp. 822--831 (1985).

\bibitem{AD}
P.~C.~Argyres and M. R. Douglas,
``New phenomena in SU(3) supersymmetric gauge theory,'' Nucl. Phys. B 448.1-2 (1995): 93-126.

\bibitem{Aspinwall:2001dz} 
  P.~S.~Aspinwall and M.~R.~Douglas,
  ``D-brane stability and monodromy,''
  JHEP {\bf 0205}, 031 (2002)
  [hep-th/0110071].

\bibitem{VB}
V.~Balasubramanian, 
``Statistical inference, Occam's razor, and statistical mechanics on the space of probability distributions,''
Neural computation 9.2 (1997): 349-368.

\bibitem{BSZ}
P.~Bleher, B.~Shiffman and S.~Zelditch,
``Universality and scaling of correlations between zeros on complex manifolds,''
Inv. Math. 142.2 (2000): 351-395.

\bibitem{Bousso:2000xa} 
  R.~Bousso and J.~Polchinski,
  ``Quantization of four form fluxes and dynamical neutralization of the cosmological constant,''
  JHEP {\bf 0006}, 006 (2000)
  [hep-th/0004134].

\bibitem{Bridgeland}
T. Bridgeland, 
``Stability conditions on triangulated categories,'' Ann. Math. (2007): 317-345.

\bibitem{Clark}
A.~Clark, 
``Whatever next? Predictive brains, situated agents, and the future of cognitive science,''
Behavioral and brain sciences 36.3 (2013): 181-204.

\bibitem{Connes:1997cr} 
  A.~Connes, M.~R.~Douglas and A.~S.~Schwarz,
  ``Noncommutative geometry and matrix theory: Compactification on tori,''
  JHEP {\bf 9802}, 003 (1998)
  [hep-th/9711162].

\bibitem{Clay} {\it Dirichlet Branes and Mirror Symmetry,}
   P. Aspinwall, M. R. Douglas, M. Gross {\it et al},
  AMS 2009.

\bibitem{double}
M.~R.~Douglas and S.~H.~Shenker, Nuclear Physics B 335.3 (1990): 635-654;
E.~Brezin and V. A. Kazakov. ``Exactly solvable field theories of closed strings'';
D.~J.~Gross and A. A. Migdal,
``Nonperturbative two-dimensional quantum gravity,'' Phys. Rev. Letters 64.2 (1990): 127.

\bibitem{Douglas:1989dd} 
  M.~R.~Douglas,
  ``Strings in Less Than One-dimension and the Generalized $KdV$ Hierarchies,''
  Phys.\ Lett.\ B {\bf 238}, 176 (1990).

\bibitem{Douglas:1996sw} 
  M.~R.~Douglas and G.~W.~Moore,
  ``D-branes, quivers, and ALE instantons,''
  hep-th/9603167.

\bibitem{Douglas:1996vj} 
  M.~R.~Douglas,
  ``Superstring dualities, Dirichlet branes and the small scale structure of space,''
  in {\it Les Houches 1995, Quantum symmetries,} eds.
A. Connes and K. Gaw\c{e}dzki, pp. 519--543, North Holland, 1998,
  hep-th/9610041.

\bibitem{Douglas:2000ah} 
  M.~R.~Douglas, B.~Fiol and C.~Romelsberger,
  ``Stability and BPS branes,''
  JHEP {\bf 0509}, 006 (2005)
  [hep-th/0002037].

\bibitem{Douglas:2000gi} 
  M.~R.~Douglas,
  ``D-branes, categories and N=1 supersymmetry,''
  J.\ Math.\ Phys.\  {\bf 42}, 2818 (2001)
  [hep-th/0011017].

\bibitem{Douglas:2002fj} 
  M.~R.~Douglas,
  ``Dirichlet branes, homological mirror symmetry, and stability,''
  Proceedings of the ICM 2002 Vol 3 pp. 395--408, Higher Education Press, Beijing, 2002,
  [math/0207021 [math-ag]].

\bibitem{Douglas:2003um} 
  M.~R.~Douglas,
  ``The Statistics of string / M theory vacua,''
  JHEP {\bf 0305}, 046 (2003)
  [hep-th/0303194].

\bibitem{DSZ}
M.~R.~Douglas, B.~Shiffman and S.~Zelditch,
``Critical points and supersymmetric vacua I,''
 Comm. Math. Phys. 252.1-3 (2004): 325-358,
 [math/0402326].

\bibitem{Douglas:2006es} 
  M.~R.~Douglas and S.~Kachru,
  ``Flux compactification,''
  Rev.\ Mod.\ Phys.\  {\bf 79}, 733 (2007)
  [hep-th/0610102].

\bibitem{Douglas:2012bu} 
  M.~R.~Douglas,
  ``The String landscape and low energy supersymmetry,''
  arXiv:1204.6626 [hep-th].

\bibitem{Douglas:2015aga} 
  M.~R.~Douglas,
  ``Calabi–Yau metrics and string compactification,''
  Nucl.\ Phys.\ B {\bf 898}, 667 (2015)
  [arXiv:1503.02899 [hep-th]].

\bibitem{Freed} {\it Five Lectures on Supersymmetry}, D. Freed, AMS 2000.

\bibitem {FS} D.~Friedan and S.~H.~Shenker,
``The analytic geometry of two-dimensional conformal field theory,''
Nuclear Physics B 281.3-4 (1987): 509-545.

\bibitem{Hitchin:1986ea} 
  N.~J.~Hitchin, A.~Karlhede, U.~Lindstrom and M.~Rocek,
  ``Hyperkahler Metrics and Supersymmetry,''
  Commun.\ Math.\ Phys.\  {\bf 108}, 535 (1987).


\bibitem{Hopfield}
J.~J.~Hopfield, 
``Neural networks and physical systems with emergent collective computational abilities,''
PNAS 79.8 (1982): 2554-2558.

\bibitem{King}
A.~D.~King, 
``Moduli of representations of finite dimensional algebras,'' 
The Quarterly Journal of Mathematics 45.4 (1994): 515-530.

\bibitem{K}
M.~Kontsevich, 
``Homological algebra of mirror symmetry,''
Proceedings of the International congress of mathematicians. Birkh\"auser, Basel, 1995.

\bibitem{LM} T.~S.~Lee and David Mumford,
``Hierarchical Bayesian inference in the visual cortex,'' JOSA A 20.7 (2003): 1434-1448.

\bibitem{MacKay}
{\it Information theory, inference and learning algorithms},
David~J.~C.~MacKay, Cambridge University Press, 2003.

\bibitem{MM} {\it Information, physics, and computation},
M.~Mezard and A. Montanari.  Oxford University Press, 2009.

\bibitem{MC} David Mumford, {\it Curves and their Jacobians}, Univ. of Michigan Press, 1975.

\bibitem{ME} David Mumford, ``Towards an enumerative geometry of the moduli space of curves''
{\it Arithmetic and geometry}, Birkh\"auser, Boston, MA, 1983. 271-328.

\bibitem {MT} David Mumford, 
``Tata lectures on theta I,'' Birkh\"auser, Basel-Boston, 1983.

\bibitem{MN} David Mumford, 
``On the computational architecture of the neocortex.'' Biological cybernetics 65.2 (1991): 135-145.

\bibitem{MP}
David Mumford, 
``Pattern theory: the mathematics of perception,'' 
Proceedings of the ICM 2002 Vol I pp. 401-422, Higher Education Press, Beijing, 2002,
 [math/0212400].

\bibitem{MP2}
{\it Pattern theory: the stochastic analysis of real-world signals},
David Mumford  and Agn\`es Desolneux,
 AK Peters/CRC Press, 2010.
 
\bibitem{Pol}
J.~Polchinski, 
``Dirichlet branes and Ramond-Ramond charges,''
 Phys. Rev. Lett. 75.26 (1995): 4724.

\bibitem{Polchinski:2006gy} 
  J.~Polchinski,
  ``The Cosmological Constant and the String Landscape,''
  hep-th/0603249.

\bibitem{Sharpe}
E.~Sharpe, ``Kahler cone substructure,'' hep-th/9810064 (1998).

\bibitem{SW} G.~J.~Sussman and J.~Wisdom. 
``Chaotic evolution of the solar system'' Science 257.5066 (1992): 56-62.

\bibitem{Thomas:1999ic} 
  R.~P.~Thomas,
  ``Derived categories for the working mathematician,''
  AMS/IP Stud.\ Adv.\ Math.\  {\bf 23}, 349 (2001)
  [math/0001045 [math-ag]].

\bibitem{Weinberg}
S.~Weinberg, 
``The cosmological constant problem,''
Rev. Mod. Phys. 61.1 (1989): 1.

\bibitem{Witten:1990hr} 
  E.~Witten,
  ``Two-dimensional gravity and intersection theory on moduli space,''
  Surveys Diff.\ Geom.\  {\bf 1}, 243 (1991).
  doi:10.4310/SDG.1990.v1.n1.a5

\end{thebibliography}

\end{document}